\documentclass[12pt,a4paper]{article}
\usepackage[utf8]{inputenc}
\usepackage{amsmath}
\usepackage{amsfonts}
\usepackage{amssymb}
\usepackage{bm}
\usepackage{float}
\usepackage{longtable}
\usepackage{pdflscape}
\usepackage[margin=1in,footskip=0.25in]{geometry}
\usepackage{graphicx}
\usepackage{multirow}
\usepackage[table,xcdraw]{xcolor}

\begin{document}
\noindent \textbf{\Large{Further study on inferential aspects of log-Lindley distribution with an application of stress-strength reliability in insurance}}\\
\\
Aniket Biswas\textsuperscript{a} \,\, ,\,\, Subrata Chakraborty\textsuperscript{b}\,\, and\,\,Meghna Mukherjee\textsuperscript{c}\\
\textsuperscript{a,b}Department of Statistics, Dibrugarh University, Dibrugarh-786004, Assam, India.
\\
\textsuperscript{c} Department of Community Medicine, IPGME\&R \& SSKM Hospitals, Kolkata-700020, India.
\\
Email:\,\,\textsuperscript{a}\textit{biswasaniket44@gmail.com} \,\,\,\,\,\textsuperscript{b}\textit{subrata\_stats@dibru.ac.in}\,\,\, \textsuperscript{c}\textit{meghnastat@gmail.com}

\section*{Abstract} The log-Lindley distribution was recently introduced in the literature as a viable alternative to the Beta distribution. This distribution has a simple structure  and possesses useful theoretical properties relevant in insurance. Classical estimation methods have been well studied. We introduce estimation of parameters from Bayesian point of view for this distribution. Explicit structure of stress-strength reliability and its inference under both classical and Bayesian set-up is addressed. Extensive simulation studies show marked improvement with Bayesian approach over classical given reasonable prior information. An application of a useful metric of discrepancy derived from stress-strength reliability is considered and computed for two categories of firm with respect to a certain financial indicator.
\\
\\
\textbf{Keywords}: log-Lindley, stress-strength reliability, risk management, Bayesian estimation, captive insurance, exponential integral. 

\section{Introduction} 
Relevance of distribution with bounded support especially in unit interval is well established for modelling proportions, rates, scores and indices etc. A number of distributions have been introduced and applied in different areas like econometrics, finance, public health among others. The log-Lindley distribution, recently introduced by Gomez-deniz et al.(2014) is simple and enjoys many properties useful in insurance and inventory management applications. Jodra and Jimenez-Gamero(2016) pointed out a crucial issue with the maximum likelihood estimation for the parameters of log-Lindley distribution and presented the following re-parametrized structure of the density function to resolve the same:
\begin{equation}
f(x|\sigma,\pi)=\sigma[\pi+\sigma(\pi-1)\log(x)]x^{\sigma-1},\quad\quad 0<x<1,\,\,\sigma>0,\,\, 0\leq\pi\leq1
\end{equation}
From Table 1 and Table 3 in Jodra and Jimenez-Gamero(2016), ML estimates for (1.1) outperforms the same for the earlier introduced distribution in terms of MSE and bias. Moreover, providing close starting values for the parameters in optimization routine failed to yield acceptable estimates for the original version as mentioned in Jodra and Jimenez-Gamero(2016).  
In this article, we work with the above re-parametrized version and refer the same as $LL(\sigma,\pi)$. If $X\sim LL(\sigma,\pi)$, then distribution function of $X$ is given by:
\begin{equation}
F(x|\sigma,\pi)=[1+\sigma(\pi-1)\log(x)]x^{\sigma},\quad\quad 0<x<1, \,\, \sigma>0,\,\, 0\leq \pi\leq 1
\end{equation}
Closed form expression of quantile function involving Lambert's $W$-function (see, Corless et al., 1996) and moments are available. Exact solution of likelihood equations are inconvenient. However, numerical optimization routine for ML estimates is available and provides satisfactory result. To the best of our knowledge, no attempts have been made so far to obtain the Bayes' estimate for the parameters when adequate prior information is available. In this article, a closed form Bayes' estimate is reported which performs very well in varied situations.\\
\\
\indent Stress-strength reliability is a stochastic measure generally applied in the the fields of engineering, medicine, econometrics, public health and allied branches. Let $X$ be strength of a system subjected to a random stress $Y$. For $X\sim LL(\sigma_1,\pi_1)$ and $Y\sim LL(\sigma_2,\pi_2)$, assuming independence between $X$ and $Y$, stress-strength reliability is defined as:
\begin{equation*}
\mathcal{R}=\int_0^{1} f_X(x|\sigma_1,\pi_1) F_Y(x|\sigma_2,\pi_2)\,dx
\end{equation*}  
A large amount of literature on stress-strength reliability estimation is available. Readers may see Krishna et al.(2017) for details. For distributions with support in unit interval, only a handful of notable works that we encounter are Kumaraswamy distribution (see Nadar et al., 2014) and Topp-Leone distribution (see Genc, 2013). Recently, Biswas and Chakraborty (2019) studied inferential issues of stress-strength reliability for unit-Lindley distribution (see, Mazucheli et al., 2019) and presented an offbeat application in public health domain. For log-Lindley distribution, expression of $\mathcal{R}$ is derived here and its inferential issues are studied in detail from classical as well as Bayesian perspective. \\
\\
\indent Gomez-Deniz et al. (2014) and Jodra and Jimenez-Gamero (2016) analysed the risk survey data of Schmit and Roth (1990) to show that, the response variable is better modelled by log-Lindley as compared to Beta, considering the same data-set with and without covariates. Here, we provide an innovative application of stress-strength reliability to measure discrepancy between two naturally existing groups of organizations induced by owning captive insurance company with respect to cost-effectiveness of risk management. Consider two groups: $A$ comprising of $n$ entities and $B$ with $m$ entities on which iid measurements on a certain indicator is available. Biswas and Chakraborty (2019) presented a discrepancy measure between $A$ and $B$ as 
\begin{equation*}
\mathcal{D}(A,B)=\mathcal{R}(B,A)-\mathcal{R}(A,B)
\end{equation*}
where, $\mathcal{R}(A,B)=P(Y<X)$. For properties and interpretation of this measure see Biswas and Chakraborty (2019).\\
\\
\indent Rest of the article is organized as follows: In section 2, we provide a short recap on maximum likelihood estimation of parameters $\sigma$ and $\pi$ though the principal focus will be on Bayes' estimation of the parameters along with computational routine for the same. In section 3, we derive the exact expression of $\mathcal{R}$ and address corresponding inferential issues. In the next section, we evaluate performance of the estimators through simulation study and interpret the findings. In section 5, an application of the metric $\mathcal{D}$ is furnished by considering the cost effectiveness of risk management data-set considered by Gomez-Deniz et al. (2014). In the last section, we provide a discussion with future direction. An appendix of relevant tables and figures is also provided.

\section{Estimation of $(\sigma,\pi)$}
Here we first discuss maximum likelihood estimation of $(\sigma,\pi)$ following Jodra and Jimenez-Gamero (2016) and additionally we derive the asymptotic covariance matrix for the estimators. Then, we derive posterior densities of $\sigma$ and $\pi$ as finite mixtures of Gamma and Beta densities, respectively. Closed form analytical expressions of Bayes' estimators of both the parameters are obtained. Issues related to calculation of exact estimates and use of statistical software is also outlined.
\subsection{Maximum likelihood estimation}
Let $\bm{X}=(X_1, X_2,..., X_m)$ be a random sample of size $m$ from $LL(\sigma,\pi)$ with pdf given in (1.1). Then, for observed $\bm{X}:=\bm{x}=(x_1, x_2, ..., x_m)$, the likelihood and log-likelihood functions of $(\sigma,\pi)$ are given respectively by:
\begin{equation}
L(\sigma,\pi|\bm{x})=\sigma^m\,\prod_{i=1}^m[\pi+\sigma(\pi-1)\log x_i]\,(\prod_{i=1}^m x_i)^{\sigma-1}
\end{equation}
\begin{equation}
\textrm{and}\quad\quad\mathcal{L}(\sigma,\pi)= m\log\sigma+(\sigma-1)\sum_{i=1}^m\log x_i +\sum_{i=1}^m \log[\pi+\sigma(\pi-1)\log x_i]
\end{equation} 
The corresponding score functions are:
\begin{eqnarray}
\frac{\delta}{\delta \sigma}\mathcal{L}(\sigma,\pi)&=&\frac{m}{\sigma}+\sum_{i=1}^m\log x_i +\sum_{i=1}^m \frac{(\pi-1)\log x_i}{\pi+\sigma(\pi-1)\log x_i}\\
\frac{\delta}{\delta \pi}\mathcal{L}(\sigma,\pi)&=&\sum_{i=1}^m\frac{1+\sigma\log x_i}{\pi+\sigma(\pi-1)\log x_i}
\end{eqnarray}
Exact maximum likelihood estimates cannot be obtained solving the above maximum likelihood equations. Thus, log-likelihood function given in (2.4) needs to be numerically maximized with respect to $\sigma$ and $\pi$ to get the estimates $\hat{\sigma}_M$ and $\hat{\pi}_M$, respectively. Jodra and Jimenez-Gamero (2016) provides data-driven starting values of the parameters for carrying out available optimization routines. For different combinations of $(\sigma,\pi)$ with only three large choices of $m$, performance of the estimator is reported therein. Therefore, it is important to investigate further for assessing performance of the same estimator for small and moderate sample sizes. Additionally, asymptotic confidence interval is reported for both the parameters. Now from (2.5) and (2.6),
\begin{eqnarray}
\frac{\delta^2}{\delta\sigma\delta\pi}\mathcal{L}(\sigma,\pi)&=&\sum_{i=1}^m\frac{\log x_i}{[\pi+\sigma(\pi-1)\log x_i]^2}\\
\frac{\delta^2}{\delta\sigma^2}\mathcal{L}(\sigma,\pi)&=& -\left[\frac{m}{\sigma^2}+\sum_{i=1}^m\left[\frac{(\pi-1)\log x_i}{\pi+\sigma(\pi-1)\log x_i}\right]^2\right]\\
\frac{\delta^2}{\delta\pi^2}\mathcal{L}(\sigma,\pi)&=&-\sum_{i=1}^m\left[\frac{1+\sigma\log x_i}{\pi+\sigma(\pi-1)\log x_i}\right]^2
\end{eqnarray}
Construct a matrix $J=((J_{ij}))$ of order $2\times 2$, where $J_{11}$, $J_{12}=J_{21}$ and $J_{22}$ are respectively given in (2.8), (2.7) and (2.9). Hence, the Fisher's information matrix is $I=((I_{ij}))_{2\times2}$,   where $I_{ij}=-E(J_{ij})$. Derivation of $E(J_{ij})$ involves the integral:
\begin{equation*}
\int_{0}^{\infty} \frac{v^k\,\exp(-\sigma v)}{\pi+\sigma(1-\pi)v}\,\,dv\quad\quad \textrm{for}\quad\quad \sigma>0,\,0<\pi<1,\, k=0,1,2.
\end{equation*}
This can be given in a closed form expression $g$ by using two popular special functions viz. Exponential Integral and Gamma function (see Abramowitz and Stegun, 1965) as follows :
\begin{equation*}
g(\sigma,\pi,k)=\exp\left(\frac{\pi}{1-\pi}\right)\,\sigma^{-(k+1)}\, EI\left(1+k,\frac{\pi}{1-\pi}\right)\, \Gamma(1+k)
\end{equation*}
where, 
\begin{eqnarray*}
EI(m,z)&=&\int_1^{\infty}\frac{\exp(-tz)}{t^m}\,dt\\
\textrm{and}\quad\quad\Gamma(z)&=&\int_{0}^{\infty} t^{z-1}\,\exp(-t)\,dt.
\end{eqnarray*}
Using $g$, $E(J_{11})=-\frac{m}{\sigma^2}-\sigma(1-\pi)^2\,g(\sigma,\pi,2)$, $E(J_{12})=-m\sigma\, g(\sigma,\pi,1)=E(J_{21})$ and $E(J_{22})=-m\sigma\,\left[ g(\sigma,\pi,0)-2\sigma\,g(\sigma,\pi,1)+\sigma^2\,g(\sigma,\pi,2)\right]$. Thus, $I$ can easily be found, from which we get covariance matrix of the vector of estimators $(\hat{\sigma}_M,\hat{\pi}_M)^\prime$, $\Sigma
=\Sigma_m(\sigma,\pi)$ by inverting $I$. Now, we provide the exact expression of asymptotic covariance and variances:
 \begin{eqnarray}
 \Sigma_{11}&=& var_m(\hat{\sigma}_M)= I_{22}\, /G(\sigma,\pi)\\
 \Sigma_{12}&=&\Sigma_{21}=cov_m(\hat{\sigma}_M,\hat{\pi}_M)= -I_{12}\, /G(\sigma,\pi)\\
 \textrm{and}\quad\quad\Sigma_{22}&=& var_m(\hat{\pi}_M)= I_{11}\, /G(\sigma,\pi)\, ,
 \end{eqnarray}
where,
\begin{equation*}
G(\sigma,\pi)=I_{11}(\sigma,\pi)\,I_{22}(\sigma,\pi)-[I_{12}(\sigma,\pi)]^2.
\end{equation*}
Under regularity conditions, the joint asymptotic distribution of the estimators:
\begin{equation}
(\hat{\sigma}_M\,,\,\hat{\pi}_M)^\prime\sim N_2((\sigma\,,\,\pi)^\prime,\Sigma_m(\sigma,\pi))
\end{equation}
Evidently, the asymptotic marginal distributions of both $\hat{\sigma}_M$ and $\hat{\pi}_M$ are both Gaussian. Based on this fact, marginal $100(1-\alpha)\%$ confidence intervals (CI) can be constructed for $\sigma$ and $\pi$ as given below:
\begin{eqnarray}
CI(\sigma)&=&\left[\hat{\sigma}_M-Z_{\alpha/2}\sqrt{\hat{var}_m(\hat{\sigma}_M)}\,\, ,\,\, \hat{\sigma}_M+Z_{\alpha/2}\sqrt{\hat{var}_m(\hat{\sigma}_M)}\right]\\
CI(\pi)&=&\left[\hat{\pi}_M-Z_{\alpha/2}\sqrt{\hat{var}_m(\hat{\pi}_M)}\,\, ,\,\, \hat{\pi}_M+Z_{\alpha/2}\sqrt{\hat{var}_m(\hat{\pi}_M)}\right]
\end{eqnarray}
where, $Z_\alpha$ denotes $(1-\alpha)$-th quantile of standard normal variate. By the invariance property, $\widehat{\Sigma}_m(\sigma,\pi)=\Sigma_m(\hat{\sigma}_M,\hat{\pi}_M)$. In view of the above desirable properties viz. invariance, consistency, asymptotic normality, maximum likelihood estimates are convenient for our subsequent inferential works on $\mathcal{R}$. 
\subsection{Bayes' estimation}
Here we first provide an equivalent expression of the likelihood function in (2.3) convenient for derivation of closed form marginal posterior densities. The following result is obtained by expanding the product in (2.3) as a weighted sum of $m+1$ different statistics using algebraic patchwork. \\
\\
\textbf{Result 1.} $L(\sigma,\pi|\bm{x})=e^{-(\sigma-1)V_1}\sum_{i=0}^m \pi^{m-i}(1-\pi)^i \sigma^{m+i} V_i$\\
where,
\begin{equation*}
V_i=
\begin{cases}
1,\quad\quad &i=0\\
\\
(-1)^i\sum_{k_1=1}^m\sum_{k_2=k_1+1}^m...\sum_{k_i=k_{i-1}+1}^m \left[\prod_{j=1}^i \log(x_{k_j})\right], \quad\quad &i=1, 2,..., (m-1)\\
\\
(-1)^m  \prod_{k=1}^m \log(x_k), \quad\quad &i=m\hspace{2cm}\blacksquare
\end{cases}
\end{equation*}
Now, we take up selection of prior on both $\sigma$ and $\pi$. As the parameter space of $\sigma$ is $(0,\infty)$, a popular choice for expressing subjective prior belief through stochastic model would be Gamma distribution. Similarly, subjective belief on $\pi$ having parameter space $(0,1)$ can be elaborated through Beta distribution. Thus we consider following independent prior distributions on the parameters.
\begin{eqnarray}
p(\sigma|\tau,\delta)&=&\frac{\tau^\delta}{\Gamma\delta} \sigma^{\delta-1} \exp(-\sigma\tau),\quad\tau\geq 0\quad\textrm{and}\quad\delta>0\\
p(\pi|\alpha,\beta)&=&\frac{\Gamma(\alpha+\beta)}{\Gamma\alpha\,\Gamma\beta}\,\pi^{\alpha-1}\,(1-\pi)^{\beta-1},\quad\alpha>0\quad\textrm{and}\quad\beta>0
\end{eqnarray}  
The joint posterior density of $(\sigma,\pi)$ is 
\begin{equation*}
p(\sigma,\pi|\bm{x})\propto L(\sigma,\pi|\bm{x})\,p(\sigma|\tau,\delta)\,p(\pi|\alpha,\beta).
\end{equation*}
Now using Result 1 along with (2.16) and (2.17), the following form is obtained
\begin{equation}
p(\sigma,\pi|\bm{x})\propto \sum_{i=0}^m V_i \left[\pi^{\alpha+m-i-1}(1-\pi)^{\beta+i-1}\right]\left[e^{-\sigma(\tau+V_1)}\sigma^{m+\tau+i-1}\right].
\end{equation}
From the structure in (2.18), it is obvious that, only some simple algebraic manipulations will turn the expressions in brackets to proper densities. We present the posterior marginals in the following result.\\
\\
\textbf{Result 2.} The marginal posterior densities of $\sigma$ and $\pi$ are respectively given by
\begin{eqnarray*}
p(\sigma|\bm{x})&=&\int_0^1 p(\sigma,\pi|\bm{x})\,d\pi=\sum_{i=0}^m w_i\,q(\sigma|\tau,\delta,\bm{x})\\
p(\pi|\bm{x})&=&\int_0^\infty p(\sigma,\pi|\bm{x})\,d\sigma=\sum_{i=0}^m w_i\,q(\pi|\alpha,\beta,\bm{x}) 
\end{eqnarray*}
where,
\begin{eqnarray*}
&&q(\sigma|\tau,\delta,\bm{x})=\frac{(\tau+V_1)^{m+\delta+i}}{\Gamma(m+\delta+i)} \sigma^{m+\delta+i-1} \exp(-\sigma(\tau+V_1))\\
\\
&&q(\pi|\alpha,\beta,\bm{x})=\frac{\Gamma(\alpha+\beta+m)}{\Gamma(\alpha+m-i)\Gamma(\beta+i)} \pi^{\alpha+m-i-1}(1-\pi)^{\beta}\\
\\
&&w_i=W_i/\sum_{i=0}^m W_i\,,\quad\quad \textrm{for}\quad\quad i=0, 1, ..., m.\\
\\
&&W_i=V_i\,\,\frac{\Gamma(\alpha+m-i)\Gamma(\beta+i)}{\Gamma(\alpha+\beta+m)}\,\,\frac{\Gamma(m+\delta+i)}{(\tau+V_1)^{m+\delta+i}}\,,\quad \textrm{for}\quad i=0, 1, ..., m.\hspace{3cm}\blacksquare 
\end{eqnarray*}
Posterior density of $\sigma$ is $(m+1)$-mixture of Gamma densities with different parameters while that of $\pi$ is $(m+1)$-mixture of different Beta densities since data and hyper-parameter dependent mixing proportions have unit sum. Posterior mean is the Bayes' estimator under squared error loss function. Hence from Result 2, we obtain the Bayes' estimators of $\sigma$ and $\pi$ respectively as
\begin{eqnarray}
\hat{\sigma}_B&=&E(\sigma|\bm{x})=\sum_{i=0}^m w_i\,\frac{m+\delta+i}{\tau+V_1}\\
\textrm{and}\quad\hat{\pi}_B&=&E(\pi|\bm{x})=\sum_{i=0}^m w_i\, \frac{\alpha+m-i}{\alpha+\beta+m}.
\end{eqnarray}
\subsection{Computational issues}
As mention in subsection 2.1, log-likelihood needs to be maximized numerically for obtaining maximum likelihood estimates. In numerical studies we used $optim$ function in \textbf{R} software with appropriate constrains.\\

\indent Under Bayesian approach, computation of $\hat{\sigma}_B$ and $\hat{\pi}_B$ is of prime interest as long as we work with squared error loss. For other loss functions, estimators will obviously change but deducing such estimator is very trivial in our case as the posterior densities are mixture of Gamma and Beta distributions. However, numerical computation of the innocent looking estimators do involve computation of $V_i$'s through the $w_i$'s. $V_i$'s involve computation of ordered sum of real numbers (see Result 1) and number of terms present in summation, $m\choose i$ explodes with sample size $m$ making exact computation of estimators beyond moderate choices of $m$, rendering it impractical. We have done numerical studies for $m$ upto $20$ and same algorithm can be followed for higher values of $m$ but we refrain from doing so due to the reason stated above.\\

\indent However, approximate Bayesian computation through STAN is practical in this kind of situation since we are sure about properness, stability and stationarity of the posteriors from Result 2. STAN applies Hamiltonian Markov chain for sampling from posterior densities and provides diagnostic measure of convergence and other details regarding the approximate Bayesian computation performed. We use STAN for drawing samples of size 10000 (after burn in 10000) for calculation of posterior means and other diagnostic measures. In simulation study, results are obtained using STAN and it is quite satisfying that the approximate computations are quite in line with the exact ones verified for small choices of $m$. Numerical and visual diagnostic outputs from STAN will be illustrated in the course real data analysis in section 5.
\section{Estimation of $\mathcal{R}$}
Let us consider independent random samples $\bm{X}= (X_1, X_2, ..., X_m)$ from $LL(\sigma_1,\pi_1)$ and $\bm{Y}= (Y_1, Y_2, ..., Y_n)$ from $LL(\sigma_2,\pi_2)$ where $m\geq n$ such that, 
\begin{equation}
\frac{n}{m}=r\in (0,1]\quad\quad\textrm{for} \quad \quad m, n \rightarrow\infty
\end{equation}
Here the stress-strength reliability parameter can be expressed as,
\begin{eqnarray}
\mathcal{R}&=&\frac{\sigma_1}{(\sigma_1+\sigma_2)^3}\left[\sigma_1(\sigma_1+3\sigma_2)+2\pi_1\sigma_2^2-2\pi_2\sigma_1\sigma_2+\pi_1\pi_2(\sigma_1\sigma_2-\sigma_2^2)\right] \nonumber \\ 
&=& f(\sigma_1,\pi_1,\sigma_2,\pi_2),\quad \textrm{say}
\end{eqnarray}
The denominator of (3.22) is non-zero for any fixed but unknown choice of $(\sigma_1,\pi_1,\sigma_2,\pi_2)$ and $f(\sigma_1,\pi_1,\sigma_2,\pi_2)$               is continuous over $(0,1)$. In the following subsection, we present maximum likelihood estimator for $\mathcal{R}$ with asymptotic confidence interval for the same. Then we approach the estimation problem from Bayesian perspective.
\subsection{Maximum likelihood estimation} 
Applying functional invariance property, the maximum likelihood estimator for $\mathcal{R}$ is given by,
\begin{equation}
\hat{\mathcal{R}}_M=f(\hat{\sigma}_{1M}, \hat{\pi}_{1M}, \hat{\sigma}_{2M}, \hat{\pi}_{1M})
\end{equation}
As in (2.13), the following asymptotic results are obvious:
\begin{eqnarray}
(\hat{\sigma}_{1M}, \hat{\pi}_{1M})^\prime&\sim& N_2((\sigma_1, \pi_1)^\prime,\Sigma_m(\sigma_1,\pi_1)\\
(\hat{\sigma}_{2M}, \hat{\pi}_{2M})^\prime &\sim& N_2((\sigma_2, \pi_2)^\prime ,\Sigma_n(\sigma_2,\pi_2)
\end{eqnarray}
\\
\textbf{Result 3} Under the assumption in (3.21), $\hat{\mathcal{R}}_M\sim N(\mathcal{R},\sigma_*^2),\quad \sigma_*^2=\sigma_{m,n}^2(\sigma_1,\pi_1,\sigma_2,\pi_2)= h_1+h_2+h_3+h_4$, where
\begin{eqnarray*}
h_1&=&a\,\frac{\sigma_1}{(\sigma_1+\sigma_2)^4}\,\left[b\,cov_n(\hat{\sigma}_{2M},\hat{\pi}_{2M})+c\,var_n(\hat{\sigma}_{2M})\right]\\
h_2&=&\frac{\sigma_1\sigma_2}{(\sigma_1+\sigma_2)^3}\,\left[(-2+\pi_1)\sigma_1-\pi_1\sigma_2\right]\,\left[b\,var_n(\hat{\pi}_{2M})+ c\,cov_n(\hat{\sigma}_{2M},\hat{\pi}_{2M}) \right]\\
h_3&=&a\,\frac{\sigma_2}{(\sigma_1+\sigma_2)^4}\,\left[d\,cov_m(\hat{\sigma}_{1M},\hat{\pi}_{1M})-e\,var_m(\hat{\sigma}_{1M})\right]\\
h_4&=&\frac{\sigma_1\sigma_2}{(\sigma_1+\sigma_2)^3}\,\left[-(-2+\pi_2)\sigma_2 +\pi_2\sigma_1\right]\,\left[d\,var_m(\hat{\pi}_{1M})- e\,cov_m(\hat{\sigma}_{1M},\hat{\pi}_{1M}) \right]\\
\end{eqnarray*}
wherein,
\begin{eqnarray*}
a&=&((-2+\pi_1)\pi2\sigma_1^2+2(-3-2\pi_1(-1+\pi_2)+2\pi_2)\sigma_1\sigma_2+\pi_1(-2+\pi_2)\sigma_2^2)\\
b&=&\frac{\sigma_1\sigma_2((-2+\pi_1)\sigma_1-\pi_1\sigma_2)}{(\sigma_1+\sigma_2)^3}\\
c&=&a\,\frac{\sigma_1}{(\sigma_1+\sigma_2)^4}\\
d&=&\frac{\sigma_1\sigma_2(\pi_2\sigma_1-(-2+\pi_2)\sigma_1)}{(\sigma_1+\sigma_2)^3}\\
e&=&a\,\frac{\sigma_2}{(\sigma_1+\sigma_2)^4}\hspace{10cm}\blacksquare
\end{eqnarray*}

Result 3, though messy, can easily be proved using delta method as follows:\\
As mentioned earlier, $f$ is continuous. Thus, applying delta method we get
\begin{equation}
f(\hat{\sigma}_{1M}, \hat{\pi}_{1M}, \hat{\sigma}_{2M}, \hat{\pi}_{1M})\sim
N((\sigma_1,\pi_1,\sigma_2,\pi_2),\sigma_*^2),\quad \textrm{asymptotically}
\end{equation}
where,
\begin{equation}
\sigma_*^2=\left[\frac{\delta f}{\delta\sigma_1}\,\,\frac{\delta f}{\delta\pi_1}\,\,\frac{\delta f}{\delta\sigma_2}\,\,\frac{\delta f}{\delta\pi_2}\right]\Sigma_{m,n}(\sigma_1,\pi_1,\sigma_2,\pi_2)\left[\frac{\delta f}{\delta\sigma_1}\,\,\frac{\delta f}{\delta\pi_1}\,\,\frac{\delta f}{\delta\sigma_2}\,\,\frac{\delta f}{\delta\pi_2}\right]^\prime
\end{equation} 
with $\frac{\delta f}{\delta\sigma_1}=-a\,\frac{\sigma_2}{(\sigma_1+\sigma_2)^4}$,
$\frac{\delta f}{\delta\pi_1}=d$,
$\frac{\delta f}{\delta\sigma_2}=a\,\frac{\sigma_1}{(\sigma_1+\sigma_2)^4}$ , 
$\frac{\delta f}{\delta\pi_2}=b$ and
\begin{center}
$$
\Sigma_{m,n}=
\begin{bmatrix}
var_m(\hat{\sigma}_{1M})& cov_m(\hat{\sigma}_{1M},\hat{\pi}_{1M}) &0&0\\
cov_m(\hat{\sigma}_{1M},\hat{\pi}_{1M}) & var_m(\hat{\pi}_{1M}) &0&0\\
0&0&var_n(\hat{\sigma}_{2M})& cov_n(\hat{\sigma}_{2M},\hat{\pi}_{2M})\\
0&0&cov_n(\hat{\sigma}_{2M},\hat{\pi}_{2M}) & var_n(\hat{\pi}_{1M})
\end{bmatrix}
$$
\end{center}
Result follows immediately by putting these in (3.27). Substituting the parameters with corresponding maximum likelihood estimates, we obtain
\begin{equation*}
\hat{\sigma}_*^2=\hat{\sigma}_{m,n}^2(\sigma_1,\pi_1,\sigma_2,\pi_2)=\sigma_{m,n}^2(\hat{\sigma}_1,\hat{\pi}_1,\hat{\sigma}_2,\hat{\pi}_2)
\end{equation*}   
Thus, using Result 3, a $100(1-\alpha)\%$ asymptotic confidence interval is given by
\begin{equation}
CI(\mathcal{R})=\left[\hat{\mathcal{R}}_M-Z_{\alpha/2}\hat{\sigma}_*\quad,\quad\hat{\mathcal{R}}_M+Z_{\alpha/2}\hat{\sigma}_* \right]
\end{equation}
\subsection{Bayes' Estimation}
Prior on $\mathcal{R}$ can sufficiently be described by prior beliefs on $(\sigma_1,\pi_1)$ and $(\sigma_2,\pi_2)$. We therefore consider same  family of prior distributions with different choices of hyper-parameters for both vector of parameters as in (2.16) and (2.17).
\begin{eqnarray*}
\sigma_1&\sim& Gamma(\tau_1,\delta_1)\,,\,\,\pi_1\sim Beta(\alpha_1,\beta_1)\hspace{0.5cm}\textrm{and}\\
\sigma_2&\sim& Gamma(\tau_2,\delta_2)\, ,\,\,\pi_2\sim Beta(\alpha_2,\beta_2).\\
\end{eqnarray*}
By plugging-in corresponding hyper-parameters in Result 2, we obtain posterior densities of $\sigma_1$, $\pi_1$, $\sigma_2$ and $\pi_2$ as $p(\sigma_1|\bm{x})$, $p(\pi_1|\bm{x})$, $p(\sigma_2|\bm{y})$ and $p(\pi_2|\bm{y})$, respectively.\\

 The computational issue raised in subsection 2.3 persists here. For moderately small sample size the mixing proportions can be computed and thus the problem of drawing sample from posterior densities reduces to that of drawing random observations from finite Gamma and Beta mixtures. In situation where at least one of the samples is large enough to the extent of inhibiting exact computation of mixing proportions, STAN is used to get four different arrays of size 10000 comprising samples from the posterior distribution of the parameters.\\
 
From the expression of $\mathcal{R}$ in (3.22), deducing the posterior distribution of $\mathcal{R}$ and hence its mean is near impossible. Therefore a sensible way-out is to draw observations from posterior distribution of $\mathcal{R}$ to compute empirical mean based on a large sample as an approximation to the true posterior mean. This is a common practice in Bayesian approach when the parametric function of interest is complicated. An algorithm for the same is presented below:\\

\textbf{Algorithm 1}
\textit{
\begin{itemize}
\item\textbf{Step 1:} Extract the $i$-th entries from 4 arrays of posterior samples viz. $\sigma_1^{(i)}$, $\pi_1^{(i)}$, $\sigma_2^{(i)}$ and $\pi_2^{(i)}$ for $i=1,  2, ..., 10000.$ 
\item\textbf{Step 2:} For each $i=1, 2, ..., 10000$, compute $\mathcal{R}^{(i)}$ using (3.22) and store it in an array.
\item\textbf{Step 3:} Calculate mean of the array obtained in Step 2 to get the posterior mean of $\mathcal{R}$, $\hat{\mathcal{R}}_B$.
\end{itemize}
}
\section{Simulation study}
As discussed in subsection 2.1, first we validate performance of maximum likelihood estimators for both the parameters through detailed simulation study. The motivation behind extending the previously done studies has already been mentioned. The experiment is carried out through the following steps:\\

\textbf{Algorithm 2}
\textit{
\begin{itemize}
\item \textbf{Step 1:} For fixed $(\sigma,\pi)$, generate a sample of size $m$ from $LL(\sigma,\pi)$ and store in array $\bm{x}$.
\item \textbf{Step 2:} For $\bm{x}$, calculate $\hat{\sigma}_M$, $\hat{\pi}_M$, $\hat{var}_m(\hat{\sigma}_M)$, $\hat{var}_m(\hat{\pi}_M)$  and thus $(\hat{\sigma}_M-\sigma)$, $(\hat{\sigma}_M-\sigma)^2$, $(\hat{\pi}_M-\pi)$, $(\hat{\pi}_M-\pi)^2$, $CI(\sigma)$ and $CI(\pi)$  and stack them into the arrays $B_\sigma$, $M_\sigma$, $B_\pi$, $M_\pi$, $C(\sigma)$ and $C(\pi)$, respectively. 
\item \textbf{Step 3:} Repeat \textbf{Step 1} to \textbf{Step 2} for $N=1000$ times.
\item \textbf{Step 4:} Average each array to get, $Bias(\hat{\sigma}_M)$, $MSE(\hat{\sigma}_M)$, $Bias(\hat{\pi}_M)$, $MSE(\hat{\pi}_M)$, $CI(\sigma)$ and $CI(\pi)$.     
\end{itemize}
}   
\noindent The results of this simulation experiment reported in Table 1 reveals the following:
\begin{itemize}
\item Bias and MSE for both the parameters reduces with increasing sample size as expected.
\item $CI(\sigma)$ shrinks with larger $m$, as desired.
\item For obvious reasons, lower limit of $CI(\pi)$ less than $0$ and upper limit of the same greater than $1$ need to be treated as $0$ and $1$, respectively. However, we reported the limits as it is to exhibit the shrinkage of $CI(\pi)$ towards the true value with growing sample size. 
\item Marginal asymptotic confidence interval for both the parameters are not performing up to the mark possibly due to the correlation between the estimators.
\item Bias in $\hat{\sigma}_M$ is increasing with $\sigma$, a quite common phenomenon for unbounded parameters.   
\end{itemize}
The results in Table 2 is produced using the following Algorithm 3.\\

\textbf{Algorithm 3}
\textit{
\begin{itemize}
\item \textbf{Step 1:} For fixed $(\sigma,\pi)$, generate a sample of size $m$ from $LL(\sigma,\pi)$ and store in array $\bm{x}$.
\item \textbf{Step 2:} Feed $(\tau_0,\delta_0,\alpha_0,\beta_0)$ along with $\bm{x}$ to STAN for generating samples from posterior distributions of both the parameters using 4 chains each with 2500 warm-up and 2500 acceptable observations on $(\sigma,\pi)$. These random samples each of size 10000 for both $\sigma$ and $\pi$ are stored in arrays $\sigma_{post}$ and $\pi_{post}$, respectively.  
\item \textbf{Step 4:} Compute 
\begin{eqnarray*}
\hat{\sigma}_B&=&mean(\sigma_{post})\\
\hat{\pi}_B&=&mean(\pi_{post})\\
L_\sigma &=& quantile(\sigma_{post},0.025)\\
U_\sigma &=& quantile(\sigma_{post},0.975)\\
L_\pi &=& quantile(\pi_{post},0.025)\\
U_\pi &=& quantile(\pi_{post},0.975)\\
\end{eqnarray*}
\item \textbf{Step 5:} Stack $(\hat{\sigma}_B-\sigma)$, $(\hat{\sigma}_B-\sigma)^2$, $(\hat{\pi}_B-\pi)$, $(\hat{\pi}_B-\pi)^2$, $(L_\sigma,U_\sigma)$ and $(L_\pi,U_\pi)$.  
\item \textbf{Step 6:} Repeat \textbf{Step 1} to \textbf{Step 5} for $N=1000$ times.
\item \textbf{Step 7:} Average each array to get, $Bias(\hat{\sigma}_B)$, $MSE(\hat{\sigma}_B)$, $Bias(\hat{\pi}_B)$, $MSE(\hat{\pi}_B)$, $Cr.I(\sigma)$ and $Cr.I(\pi)$.     
\end{itemize}
}  
\noindent Note that, for Step 2 in Algorithm 3,  reasonable choices of hyper-parameters are needed. For synthetically generated data, we take hyper-parameters in such a manner that the prior mean turns out to be the fixed value of parameter. For different values of $(\sigma,\pi)$ in our study, the hyper-parameter values are listed below:
\begin{longtable}{ccc|ccc}
$\sigma$&$\tau_0$&$\delta_0$&$\pi$&$\alpha_0$&$\beta_0$\\
\hline
1.0&1.0&1.0&0.2&1.0&4.0\\
2.5&2.0&5.0&0.5&1.0&1.0\\
3.5&2.0&7.0&0.7&3.5&1.5\\
\hline
\end{longtable}
Findings from Table 2:
\begin{itemize}
\item For small and moderate sample size ($\leq 100$), Bayes' estimators for both the parameters have strictly lower bias and MSE compared to maximum likelihood estimators.
\item Credible interval is more crisp than confidence intervals.
\item MSE of $\hat{\sigma}_B$ is decreasing with increasing $m$.
\item MSE of $\hat{\pi}_B$ increases with $m$ upto $m=75$ then decreases with increasing $m$. This may be attributed to the fact that, for smaller sample sizes, precise prior belief makes the estimator more efficient and for moderate sample sizes, the prior though perfect looses importance in the posterior inference. The decreasing trend towards larger sample sizes is mainly due to adequate information in sample to make the estimator efficient even when prior is not too precise.
\end{itemize}
Results in Table 3 are generated using Algorithm 1, 2 and 3 with necessary adjustments. Findings from Table 3:
\begin{itemize}
\item Bias and MSE of both the estimators for $\mathcal{R}$ decreases with sample sizes implying consistency.
\item It is observed that biases in $\hat{\mathcal{R}}_M$ and $\hat{\mathcal{R}}_B$ are more or less comparable.
\item $\hat{\mathcal{R}}_B$ beats $\hat{\mathcal{R}}_M$ uniformly with respect to MSE.
\item Credible intervals are much concentrated around the true value as compared to confidence intervals.
\item For fixed total sample size, maximum precision of both the estimators is attained for equal allocation (Note the cases for $(m,n)=(125,125)$ and $(m,n)=(150,100)$).       
\end{itemize} 
\section{Data analysis}
As mentioned in section 1, we consider the risk management data-set in which observations on $7$ variables are given for $73$ organizations. We do not delve into the aspects of modelling the response variable with or without covariates as the same has been adequately addressed in Gomez-Deniz et al. (2014) and Jodra and Jimenez-Gamero (2016). Out of $7$ variables, viz. FIRMCOST, CAP, ASSUME, SIZELOG, INDCOST, CENTRAL and SOPH, consider the first two variables:
\begin{itemize}
\item \textbf{FIRMCOST}: A measure of a firm's risk management cost effectiveness. 
\item \textbf{CAP}: An indicator variable. $CAP=1$ if the firm owns a captive insurance company and $CAP=0$ if it currently does not.
\end{itemize}  
There exists a natural categorization in the 73 firms which is clear when we consider $A$ to be the group comprising 25 firms owning a captive insurance company (CAP=1) and $B$ to be the same comprising 48 firms not owning a captive insurance company (CAP=0). As in the previous studies with this data, the response variable is taken to be FIRMCOST/100 for adaptation to log-Lindley model. Our main interest lies in applying the metric $\mathcal{D}(A,B)$ with related inference for standardized probabilistic quantification of the discrepancy between $A$ and $B$ with respect to cost effectiveness, an important financial indicator. In accordance with our discussion in section 3, $X$ stands for the random response for firms in group $A$ and $Y$ stands for the same for group $B$. Estimate of the relevant parameters along with variances and confidence intervals are presented below:
\begin{longtable}{cccc}
Parameter&ML Estimate& Variance & 95\%CI\\
\hline
$\sigma_1$&0.7282&0.0055&(0.5728,0.8635)\\
$\pi_1$&0.0000&0.0012&(0.0000,0.0679)\\
$\sigma_2$&0.6301&0.0127&(0.4092,0.8510)\\
$\pi_2$&0.1029&0.0444&(0.0000,0.5159)\\
\hline
$\mathcal{R}(A,B)$&0.5216&0.0125&(0.3025,0.7407)\\
$\mathcal{D}(A,B)$&-0.0432&0.0500&(-0.4815,0.3951)\\
\hline
\end{longtable}
From definition of $\mathcal{D}(A,B)$, it is obvious that $\mathcal{D}(A,B)=0$ indicates almost similar situation for the two groups. In our context of study, a negative value of the same indicates a better situation in $B$ than in $A$ as the response is defined as total property and casualty premiums and uninsured losses as a proportion of total losses, an indicator whose reduction is every manager's priority. In fact, from estimated $\mathcal{D}(A,B)$, we can say that firms owning captive insurance company is in a slightly better situation with respect to cost effectiveness of its risk management. Though the discrepancy turns out  be statistically insignificant in this case, it may become significant in other similar situations with same or different choice of financial indicator.\\

As observed in simulation studies, the Bayes' estimator out-performs the maximum likelihood estimator when the experimenter can provide a reasonable choice of hyper-parameters for depiction of the true picture. With the notations in subsection 3.2, we choose $\alpha_1=0.001$, $\beta_1=5.5$, $\tau_1=5.0$, $\delta_1=3.5$, $\alpha_2=0.5$, $\beta_2=5.5$, $\tau_2=5.0$ and $\delta_2=3.5$. Under the above setting, Bayes' estimate of the relevant parameters along with credible intervals are presented below:
\begin{longtable}{ccc}
Parameter&Bayes' Estimate & 95\%Cr.I\\
\hline
$\sigma_1$&0.7165&(0.5837,0.8632)\\
$\pi_1$&0.0014&(0.0000,0.0184)\\
$\sigma_2$&0.6395&(0.4699,0.8304)\\
$\pi_2$&0.0866&(0.0003,0.3000)\\
\hline
$\mathcal{R}(A,B)$&0.5213&(0.3985,0.6440)\\
$\mathcal{D}(A,B)$&-0.0425&(-0.2880,0.2030)\\
\hline
\end{longtable}
There is no convergence issue as indicated by $Rhat=1$, provided by STAN. Trace-plots for posterior sampling from $(\sigma_1,\pi_1)$ is given in Figure 1 and the same for $(\sigma_2,\pi_2)$ in Figure 2, given in Appendix.
\section{Discussion}
The cost effectiveness or performance risk management of a large organization depend on variety of factors both financial and non-financial. One of the major factor is the use (or otherwise) of a captive insurer by the concerned organization. A general intuition tells us that, employing a captive ensurer is expected to yield lower per unit risk management cost. Thus, the investigation of the problem of measuring the discrepancy between these categories considered in this work is well justified and useful.\\

There are some technical limitations in this paper. Exact computation of Bayes' estimator from supplied closed form expression remains a computational challenge for large sample size. Using improper prior $p(\sigma,\pi)\propto 1$, Bayes' estimators were derived but not reported due to poor performance whereas, implementing Jeffrey's prior is found to be practically inconvenient due to computational complexities. Thus, we recommend the practitioners to go the Bayesian way only when informative prior knowledge is available. Resolving these issues may pave ways for future investigations. Moreover, consideration of dependence between two groups will necessitate a fresh look into the inferential aspects of stress-strength reliability.
\section*{Reference}
Abramowitz, M., \& Stegun, I. A. (1965). \textit{Handbook of mathematical functions: with formulas, graphs, and mathematical tables (Vol. 55)}. Courier Corporation.\\
\\
Biswas, A., \& Chakraborty, S. (2019). $\mathcal{R}=P(Y<X)$ for unit-Lindley distribution: inference with an application in public health. \textit{https://arxiv.org/abs/1904.06181}\\
\\
Corless, R. M., Gonnet, G. H., Hare, D. E., Jeffrey, D. J., \& Knuth, D. E. (1996). On the LambertW function. \textit{Advances in Computational mathematics}, 5(1), 329-359.\\
\\
Genç,A.I. (2013). Estimation of $P(X>Y)$ with Topp–Leone distribution, \textit{Journal of
Statistical Computation and Simulation}, 83:2, 326-339.\\
\\
Gómez-Déniz, E., Sordo, M. A., \& Calderín-Ojeda, E. (2014). The Log–Lindley distribution as an alternative to the beta regression model with applications in insurance. \textit{Insurance: Mathematics and Economics}, 54, 49-57.\\
\\
Jodrá, P., \& Jiménez-Gamero, M. D. (2016). A note on the Log-Lindley distribution. \textit{Insurance: Mathematics and Economics}, 71, 189-194.\\
\\
Krishna, H., Dube, M., \& Garg, R. (2017). Estimation of $P(Y< X)$ for progressively first-failure-censored generalized inverted exponential distribution.\textit{Journal of Statistical Computation and Simulation}, 87(11), 2274-2289.\\
\\
Mazucheli, J., Menezes, A. F. B., \& Chakraborty, S. (2019). On the one parameter unit-Lindley distribution and its associated regression model for proportion data. \textit{Journal of Applied Statistics}, 46(4), 700-714.\\
\\
Nadar, M., Kızılaslan, F., \& Papadopoulos, A. (2014). Classical and Bayesian estimation of $P (Y<X)$ for Kumaraswamy's distribution. \textit{Journal of Statistical Computation and Simulation}, 84(7), 1505-1529.\\
\\
Schmit, J. T., \& Roth, K. (1990). Cost effectiveness of risk management practices. \textit{Journal of Risk and Insurance}, 455-470.\\
\section*{Appendix}

\begin{longtable}{ccccccc}
\caption{Bias, MSE and CI of ML estimator for $\sigma$ and$\pi$}\\
\endfirsthead
\hline
\multicolumn{7}{c}{$\sigma=1.0\quad,\quad \pi=0.2$}\\
\hline
$m$&Bias$(\hat{\sigma}_M)$&MSE$(\hat{\sigma}_M)$& CI$(\sigma)$&Bias$(\hat{\pi}_M)$&MSE$(\hat{\pi}_M)$& CI$(\pi)$\\
\hline
5&0.1600&0.2991&(0.3370,2.0228)&-0.0150&0.0921&(-0.7029,1.0086)\\
10&0.0405 &0.0979 &(0.4973,1.6330) &0.0315 &0.0980&(-0.4400,0.8067)\\
15&0.0330&0.0708&(0.5843,1.5144)&0.0137&0.0773&(-0.3476,0.7146)\\
20&0.0102&0.0545&(0.6204,1.4291)&0.0235&0.0696&(-0.2927,0.6884)\\
25&0.0203&0.0447&(0.6687,1.3948)&0.0044&0.0583&(-0.2456,0.6120)\\
50&0.0007&0.0214&(0.7400,1.2676)&0.0099&0.0331&(-0.1354,0.5457)\\
75&0.0089&0.0137&(0.7908,1.2270)&-0.0032&0.0210&(-0.0839,0.4774)\\
100&0.0032&0.0113&(0.8150,1.1961)&0.0026&0.0194&(-0.0508,0.4480)\\
125&0.0035&0.0090&(0.8350,1.1760)&-0.0032&0.0158&(-0.0293,0.4164)\\
150&0.0024&0.0078&(0.8468,1.1600)&-0.0011&0.0119&(-0.0105,0.4050)\\
\hline
\multicolumn{7}{c}{$\sigma=2.5\quad,\quad \pi=0.2$}\\
\hline
$m$&Bias$(\hat{\sigma}_M)$&MSE$(\hat{\sigma}_M)$& CI$(\sigma)$&Bias$(\hat{\pi}_M)$&MSE$(\hat{\pi}_M)$& CI$(\pi)$\\
\hline
5&0.3218&1.5206&(0.8120,4.9486)&-0.0041&0.0994&(-0.6276,0.9419)\\
10&0.1251&0.7436&(1.2457,4.1141)&0.0242&0.0907&(-0.4439,0.8119)\\
15&0.0741&0.4162&(1.4543,3.7639)&0.0122&0.0771&(-0.3295,0.6984)\\
20&0.0493&0.3307&(1.5681,3.6109)&0.0169&0.0653&(-0.2801,0.6554)\\
25&0.0532&0.2958&(1.6685,3.4986)&0.0116&0.0587&(-0.2370,0.6165)\\
50&0.0273&1.1279&(1.8721,3.1940)&0.0033&0.0335&(-0.1352,0.5322)\\
75&0.0098&0.0897&(1.9643,3.0574)&0.0048&0.0221&(-0.0819,0.4900)\\
100&0.0190&0.0629&(2.0425,2.9977)&-0.0019&0.0162&(-0.0523,0.4469)\\
125&0.0066&0.0545&(2.0806,2.9350)&-0.0022&0.0131&(-0.0280,0.4221)\\
150&0.0104&0.0452&(2.1183,2.9025)&-0.0007&0.0118&(-0.0087,0.4075)\\
\hline
\multicolumn{7}{c}{$\sigma=3.5\quad,\quad \pi=0.5$}\\
\hline
$m$&Bias$(\hat{\sigma}_M)$&MSE$(\hat{\sigma}_M)$& CI$(\sigma)$&Bias$(\hat{\pi}_M)$&MSE$(\hat{\pi}_M)$& CI$(\pi)$\\
\hline
5&1.1770&7.1055&(1.2168,8.5988)&-0.1856&0.1773&(-0.8420,1.3066)\\
10&0.5288&2.6140&(1.8128,6.6643)&-0.1092&0.1475&(-0.6925,1.2906)\\
15&0.3538&1.8329&(2.0634,6.0427)&-0.0721&0.1275&(-0.5215.1.2063)\\
20&0.3037&1.3957&(2.2507,5.7335)&-0.0618&0.1148&(-0.3687,1.0832)\\
25&0.1307&0.8987&(2.2716,5.3193)&-0.0282&0.1047&(-0.3152,1.1052)\\
50&0.0795&0.5522&(2.5690,4.7834)&-0.0067&0.0669&(-0.1015,1.0012)\\
75&0.0278&0.3811&(2.6980,4.5063)&0.0000&0.0518&(0.0121,0.9230)\\
100&-0.0119&0.3129&(2.7640,4.3295)&0.0187&0.0423&(0.0841,0.9027)\\
125&-0.0032&0.2481&(2.8312,4.2347)&0.0092&0.0334&(0.1296,0.8564)\\
150&-0.0050&0.1983&(2.8766,4.1579)&0.0052&0.0276&(0.1656,0.8245)\\
\hline
\multicolumn{7}{c}{$\sigma=1.0\quad,\quad \pi=0.7$}\\
\hline
$m$&Bias$(\hat{\sigma}_M)$&MSE$(\hat{\sigma}_M)$& CI$(\sigma)$&Bias$(\hat{\pi}_M)$&MSE$(\hat{\pi}_M)$& CI$(\pi)$\\
\hline
5&0.5305&1.0990&(0.3639,2.8512)&-0.3164&0.2577&(-0.9279,1.4895)\\
10&0.2495&0.3223&(0.5348,2.1654)&-0.1947&0.17644&(-0.7645,1.5301)\\
15&0.2172&0.2684&(0.6624,1.9871)&-0.1648&0.1604&(-0.5405,1.3426)\\
20&0.1196&0.1439&(0.6616,1.7670)&-0.1093&0.1265&(-0.4721,1.4007)\\
25&0.1073&0.1274&(0.7026,1.6994)&-0.0846&0.1101&(-0.3999,1.3893)\\
50&0.0473&0.0633&(0.7793,1.4760)&-0.0419&0.0731&(-0.0920,1.2028)\\
75&0.0281&0.0473&(0.8100,1.3668)&-0.0251&0.0568&(0.0430,1.1442)\\
100&0.0268&0.0373&(0.8384,1.3195)&-0.0161&0.0477&(0.1317,1.0964)\\
125&0.0141&0.0305&(0.8444,1.2658)&-0.0176&0.0417&(0.1917,1.0627)\\
150&0.0053&0.0272&(0.8565,1.2398)&-0.0009&0.0370&(0.2424,1.0411)\\
\hline
\multicolumn{7}{c}{$\sigma=2.5\quad,\quad \pi=0.7$}\\
\hline
$m$&Bias$(\hat{\sigma}_M)$&MSE$(\hat{\sigma}_M)$& CI$(\sigma)$&Bias$(\hat{\pi}_M)$&MSE$(\hat{\pi}_M)$& CI$(\pi)$\\
\hline
5&1.2139&5.5828&(0.9046,6.9976)&-0.3213&0.2590&(-0.9252,1.4702)\\
10&0.6797&2.3966&(1.4275,5.4695)&-0.2058&0.1938&(-0.7510,1.4577)\\
15&0.4870&1.3092&(1.5953,4.8675)&-0.1592&0.1551&(-0.5716,1.3958)\\
20&0.3783&0.9708&(1.6886,4.5067)&-0.1158&0.1265&(-0.4920,1.4300)\\
25&0.3403&0.8110&(1.7848,4.2971)&-0.1177&0.1183&(-0.3641,1.3161)\\
50&0.1142&0.3940&(1.9287,3.6555)&-0.0431&0.0717&(-0.1043,1.2213)\\
75&0.0812&0.2746&(2.0316,3.4261)&-0.0283&0.0577&(0.0408,1.1436)\\
100&0.0379&0.3326&(2.0818,3.2727)&-0.0087&0.0491&(0.1241,1.1119)\\
125&0.0207&1.1979&(2.1149,3.1701)&-0.0058&0.0436&(0.1745,1.0852)\\
150&0.0022&0.1788&(2.1372,3.0925)&-0.0038&0.0382&(0.2105,1.0794)\\
\hline
\end{longtable} 
\begin{longtable}{ccccccc}
\caption{Bias, MSE and Credible Interval of Bayes' estimator for $\sigma$ and$\pi$}\\
\endfirsthead
\hline
\multicolumn{7}{c}{$\sigma=1.0\quad,\quad \pi=0.2$}\\
\hline
$m$&Bias$(\hat{\sigma}_B)$&MSE$(\hat{\sigma}_B)$& Cr.I$(\sigma)$&Bias$(\hat{\pi}_B)$&MSE$(\hat{\pi}_B)$& Cr.I$(\pi)$\\
\hline
5&0.0991&0.1500&(0.4864,1.9176)&0.0006&0.0016&(0.0084,0.5855)\\
10&0.0751&0.0733&(0.6045,1.6457)&-0.0026&0.0029&(0.0101,0.5652)\\
15&0.0450&0.0457&(0.6495,1.5018)&-0.0029&0.0039&(0.0019,0.5517)\\
20&0.0261&0.0316&(0.6756,1.4187)&0.0021&0.0052&(0.0155,0.5446)\\
25&0.0252&0.0271&(0.7034,1.3775)&-0.0020&0.0053&(0.0159,0.5294)\\
50&0.0130&0.0152&(0.7658,1.2676)&0.0046&0.0074&(0.0278,0.4948)\\
75&0.0058&0.0106&(0.7961,1.2163)&0.0058&0.0082&(0.0376,0.4658)\\
100&0.0039&0.0084&(0.8188,1.1876)&0.0016&0.0085&(0.0442,0.4363)\\
125&0.0026&0.0068&(0.8343,1.1686)&0.0023&0.0075&(0.0510,0.4195)\\
150&0.0014&0.0059&(0.8459,1.1540)&0.0013&0.0073&(0.0561,0.4040)\\
\hline
\multicolumn{7}{c}{$\sigma=2.5\quad,\quad \pi=0.2$}\\
\hline
$m$&Bias$(\hat{\sigma}_B)$&MSE$(\hat{\sigma}_B)$& Cr.I$(\sigma)$&Bias$(\hat{\pi}_B)$&MSE$(\hat{\pi}_B)$& Cr.I$(\pi)$\\
\hline
5&0.1140&0.2989&(1.3721,4.2073)&0.0018&0.0019&(0.0088,0.5812)\\
10&0.0878&0.2288&(1.5632,3.8160)&0.0008&0.0032&(0.0108,0.5632)\\
15&0.0748&0.1899&(1.6719,3.6141)&-0.0001&0.0039&(0.0124,0.5498)\\
20&0.0613&0.1599&(1.7422,3.4794)&-0.0009&0.0049&(0.0149,0.5352)\\
25&0.0678&0.1347&(1.8054,3.4044)&-0.0016&0.0055&(0.0171,0.5233)\\
50&0.0160&0.0806&(1.9260,3.1267)&-0.0007&0.0071&(0.0261,0.4835)\\
75&0.0201&0.0596&(2.0079,3.0364)&0.0046&0.0074&(0.0375,0.4620)\\
100&0.0131&0.0469&(2.0581,2.9657)&0.0001&0.0070&(0.0429,0.4333)\\
125&0.0088&0.0393&(2.0967,2.9162)&-0.0024&0.0064&(0.0486,0.4115)\\
150&0.0102&0.0328&(2.1271,2.8872)&-0.0001&0.0064&(0.0556,0.4002)\\
\hline
\multicolumn{7}{c}{$\sigma=3.5\quad,\quad \pi=0.5$}\\
\hline
$m$&Bias$(\hat{\sigma}_B)$&MSE$(\hat{\sigma}_B)$& Cr.I$(\sigma)$&Bias$(\hat{\pi}_B)$&MSE$(\hat{\pi}_B)$& Cr.I$(\pi)$\\
\hline
5&0.1162&0.3712&(1.8741,5.8674)&0.0069&0.0100&(0.0530,0.9651)\\
10&0.0711&0.3516&(2.0286,5.4618)&0.0120&0.0147&(0.0745,0.9587)\\
15&0.0305&0.3034&(2.1119,5.2104)&0.0174&0.0172&(0.0944,0.9540)\\
20&0.0248&0.2782&(2.1795,5.0786)&0.0222&0.0178&(0.1078,0.9510)\\
25&0.0352&0.2499&(2.2425,4.9943)&0.0197&0.0197&(0.1184,0.9455)\\
50&-0.0246&0.2036&(2.3721,4.6503)&0.0331&0.0230&(0.1747,0.9256)\\
75&-0.0686&0.1865&(2.4324,4.4437)&0.0412&0.0238&(0.2115,0.9110)\\
100&-0.0769&0.1627&(2.4906,4.3317)&0.0361&0.0215&(0.2316,0.8930)\\
125&-0.0809&0.1508&(2.5386,4.2588)&0.0414&0.0214&(0.2555,0.8812)\\
150&-0.1080&0.1315&(2.5558,4.1726)&0.0476&0.0200&(0.2768,0.8751)\\
\hline
\multicolumn{7}{c}{$\sigma=1.0\quad,\quad \pi=0.7$}\\
\hline
$m$&Bias$(\hat{\sigma}_B)$&MSE$(\hat{\sigma}_B)$& Cr.I$(\sigma)$&Bias$(\hat{\pi}_B)$&MSE$(\hat{\pi}_B)$& Cr.I$(\pi)$\\
\hline
5&0.1327&0.2127&(0.4161,2.1822)&0.0017&0.0008&(0.2970,0.9711)\\
10&0.0744&0.1045&(0.5234,1.8130)&0.0067&0.0014&(0.3151,0.9711)\\
15&0.0536&0.0712&(0.5770,1.6671)&0.0008&0.0021&(0.3152,0.9699)\\
20&0.0415&0.0506&(0.6121,1.5829)&0.0052&0.0026&(0.3300,0.9700)\\
25&0.0227&0.0397&(0.6286,1.5108)&0.0032&0.0035&(0.3351,0.9692)\\
50&0.0119&0.0209&(0.6955,1.3894)&0.0823&0.0053&(0.3764,0.9680)\\
75&-0.0035&0.0150&(0.7182,1.3211)&0.0174&0.0068&(0.4048,0.9659)\\
100&-0.0053&0.0131&(0.7363,1.2893)&0.0164&0.0082&(0.4208,0.9624)\\
125&-0.0115&0.0103&(0.7463,1.2620)&0.0247&0.0070&(0.4417,0.9635)\\
150&-0.0169&0.0094&(0.7528,1.2388)&0.0257&0.0079&(0.4553,0.9614)\\
\hline
\multicolumn{7}{c}{$\sigma=2.5\quad,\quad \pi=0.7$}\\
\hline
$m$&Bias$(\hat{\sigma}_B)$&MSE$(\hat{\sigma}_B)$& Cr.I$(\sigma)$&Bias$(\hat{\pi}_B)$&MSE$(\hat{\pi}_B)$& Cr.I$(\pi)$\\
\hline
5&0.0991&0.2987&(1.2475,4.4230)&0.0056&0.0015&(0.3100,0.9703)\\
10&0.0906&0.2520&(1.4273,4.0883)&0.0043&0.0019&(0.3177,0.9696)\\
15&0.0928&0.2227&(1.5373,3.9169)&0.0051&0.0026&(0.3284,0.9691)\\
20&0.0498&0.1753&(1.5866,3.7383)&0.0082&0.0029&(0.3414,0.9687)\\
25&0.0527&0.1700&(1.6418,3.6644)&0.0065&0.0035&(0.3459,0.9687)\\
50&0.0052&0.0981&(1.7582,3.3883)&0.0120&0.0052&(0.3823,0.9662)\\
75&-0.0192&0.0757&(1.8122,3.2562)&0.0176&0.1160&(0.4087,0.9648)\\
100&-0.0281&0.0638&(1.8533,3.1807)&0.0256&0.0069&(0.4346,0.9639)\\
125&-0.0307&0.0619&(1.8800,3.1343)&0.0260&0.0068&(0.4459,0.9625)\\
150&-0.0315&0.0561&(1.9001,3.0964)&0.0206&0.0072&(0.4509,0.9583)\\
\hline
\end{longtable} 
\newpage
\begin{longtable}{ccccccc}
\caption{Bias, MSE and CI/Credible Interval of maximum likelihood/ Bayes' estimator for $\mathcal{R}$}\\
\endfirsthead
\hline
\multicolumn{7}{c}{$\sigma_1=1.0\quad,\quad \pi_1=0.2\quad,\quad\sigma_2=1.0\quad,\quad\pi_2=0.2\quad,\quad\mathcal{R}=0.5$}\\
\hline
$m,n$&Bias$(\hat{\mathcal{R}}_M)$&MSE$(\hat{\mathcal{R}}_M)$& CI$(\hat{\mathcal{R}}_M)$&Bias$(\hat{\mathcal{R}}_B)$&MSE$(\hat{\mathcal{R}}_B)$& Cr.I$(\hat{\mathcal{R}}_B)$\\
\hline
5,5&-0.0021&0.0241&(-0.0219,1.0161)&-0.0040&0.0166&(0.2322,0.7614)\\
10,5&-0.0026&0.0195&(0.0551,0.9420)&-0.0004&0.0140&(0.2677,0.7446)\\
15,10&-0.0038&0.0116&(0.1477,0.8475)&-0.0010&0.0086&(0.3111,0.6921)\\
15,15&-0.0024&0.0093&(0.1934,0.8064)&-0.0011&0.0073&(0.3253,0.6726)\\
25,25&0.0004&0.0058&(0.2457,0.7543)&-0.0023&0.0053&(0.3588,0.6368)\\
50,25&-0.0010&0.0043&(0.2741,0.7245)&-0.0007&0.0038&(0.3797,0.6234)\\
100,75&-0.0031&0.0018&(0.3559,0.6380)&-0.0004&0.0016&(0.4203,0.5793)\\
125,125&-0.0005&0.0012&(0.3857,0.6138)&-0.0016&0.0010&(0.4349,0.5683)\\
150,100&0.0002&0.0013&(0.3769,0.6236)&-0.0003&0.0011&(0.4319,0.5677)\\
150,150&0.0003&0.0010&(0.3952,0.6053)&-0.0011&0.0009&(0.4399,0.5622)\\
\hline
\multicolumn{7}{c}{$\sigma_1=1.0\quad,\quad \pi_1=0.2\quad,\quad\sigma_2=2.5\quad,\quad\pi_2=0.2\quad,\quad\mathcal{R}=0.2297$}\\
\hline
$m,n$&Bias$(\hat{\mathcal{R}}_M)$&MSE$(\hat{\mathcal{R}}_M)$& CI$(\hat{\mathcal{R}}_M)$&Bias$(\hat{\mathcal{R}}_B)$&MSE$(\hat{\mathcal{R}}_B)$& Cr.I$(\hat{\mathcal{R}}_B)$\\
\hline
5,5&0.0037&0.0152&(-0.1469,0.6081)&0.0262&0.0087&(0.0890,0.4875)\\
10,5&0.0013&0.0109&(-0.1204,0.5750)&0.0174&0.0058&(0.1028,0.4496)\\
15,10&0.0046&0.0075&(-0.0362,0.4983)&0.0122&0.0045&(0.1203,0.4011)\\
15,15&0.0073&0.0064&(-0.0054,0.4735)&0.0101&0.0039&(0.1253,0.3850)\\
25,25&0.0003&0.0038&(0.0476,0.4121)&0.0059&0.0028&(0.1409,0.3506)\\
50,25&0.0028&0.0026&(0.0691,0.3952)&0.0060&0.0018&(0.1524,0.3361)\\
100,75&0.0005&0.0012&(0.1305,0.3300)&0.0010&0.0009&(0.1728,0.2961)\\
125,125&0.0004&0.0009&(0.1467,0.3126)&0.0015&0.0008&(0.1796,0.2901)\\
150,100&0.0000&0.0008&(0.1437,0.3158)&0.0017&0.0007&(0.1813,0.2872)\\
150,150&0.0003&0.0007&(0.1545,0.3052)&0.0008&0.0006&(0.1841,0.2818)\\
\hline
\multicolumn{7}{c}{$\sigma_1=3.5\quad,\quad \pi_1=0.5\quad,\quad\sigma_2=1.0\quad,\quad\pi_2=0.7\quad,\quad\mathcal{R}=0.7576$}\\
\hline
$m,n$&Bias$(\hat{\mathcal{R}}_M)$&MSE$(\hat{\mathcal{R}}_M)$& CI$(\hat{\mathcal{R}}_M)$&Bias$(\hat{\mathcal{R}}_B)$&MSE$(\hat{\mathcal{R}}_B)$& Cr.I$(\hat{\mathcal{R}}_B)$\\
\hline
5,5&0.0102&0.0188&(0.3071,1.2381)&-0.0199&0.0072&(0.5238,0.8983)\\
10,5&0.0176&0.0157&(0.3458,1.2285)&-0.0090&0.0064&(0.5568,0.8981)\\
15,10&0.0193&0.0086&(0.4175,1.1614)&-0.0084&0.0044&(0.5990,0.8707)\\
15,15&0.0146&0.0061&(0.4342,1.1253)&-0.0076&0.0032&(0.6145,0.8600)\\
25,25&0.0088&0.0039&(0.4915,1.0590)&-0.0022&0.0023&(0.6488,0.8547)\\
50,25&0.0128&0.0032&(0.5068,1.0464)&-0.0011&0.0018&(0.6623,0.8396)\\
100,75&0.0052&0.0011&(0.6013,0.9333)&0.0005&0.0008&(0.6977,0.8146)\\
125,125&0.0034&0.0008&(0.6345,0.8957)&-0.0007&0.0005&(0.7065,0.8048)\\
150,100&0.0041&0.0009&(0.6250,0.9080)&0.0011&0.0006&(0.7069,0.8083)\\
150,150&0.0039&0.0006&(0.6466,0.8835)&0.0004&0.0005&(0.7117,0.8024)\\
\hline
\multicolumn{7}{c}{$\sigma_1=3.5\quad,\quad \pi_1=0.5\quad,\quad\sigma_2=1.0\quad,\quad\pi_2=0.7\quad,\quad\mathcal{R}=0.7576$}\\
\hline
$m,n$&Bias$(\hat{\mathcal{R}}_M)$&MSE$(\hat{\mathcal{R}}_M)$& CI$(\hat{\mathcal{R}}_M)$&Bias$(\hat{\mathcal{R}}_B)$&MSE$(\hat{\mathcal{R}}_B)$& Cr.I$(\hat{\mathcal{R}}_B)$\\
\hline
5,5&-0.0105&0.0208&(-0.268,0.7790)&0.0180&0.0079&(0.1188,0.5223)\\
10,5&-0.0186&0.0144&(-0.2027,0.7080)&0.0151&0.0064&(0.1409,0.4919)\\
15,10&-0.0112&0.0087&(-0.1347,0.6528)&0.0079&0.0043&(0.1578,0.44453)\\
15,15&-0.0103&0.0072&(-0.0893,0.6088)&0.0063&0.0046&(0.1634,0.4294)\\
25,25&-0.0108&0.0044&(-0.0195,0.5358)&0.0051&0.0026&(0.1827,0.3991)\\
50,25&-0.0062&0.0030&(-0.0103,0.5271)&0.0016&0.0019&(0.1928,0.3804)\\
100,75&-0.0031&0.0012&(0.0993,0.4431)&0.0002&0.0009&(0.2178,0.3429)\\
125,125&-0.0035&0.0008&(0.1339,0.4101)&-0.0008&0.0007&(0.2257,0.3322)\\
150,100&-0.0004&0.0009&(0.1315,0.4199)&0.0015&0.0007&(0.2272,0.3347)\\
150,150&-0.0015&0.0007&(0.1449,0.4013)&0.0005&0.0006&(0.2303,0.3283)\\
\hline
\multicolumn{7}{c}{$\sigma_1=2.5\quad,\quad \pi_1=0.7\quad,\quad\sigma_2=1.0\quad,\quad\pi_2=0.2\quad,\quad\mathcal{R}=0.8373$}\\
\hline
$m,n$&Bias$(\hat{\mathcal{R}}_M)$&MSE$(\hat{\mathcal{R}}_M)$& CI$(\hat{\mathcal{R}}_M)$&Bias$(\hat{\mathcal{R}}_B)$&MSE$(\hat{\mathcal{R}}_B)$& Cr.I$(\hat{\mathcal{R}}_B)$\\
\hline
5,5&-0.0112&0.0123&(0.4508,1.2008)&-0.0274&0.0060&(0.6088,0.9407)\\
10,5&-0.0123&0.0096&(0.4972,1.1622)&-0.0218&0.0052&(0.6394,0.9374)\\
15,10&-0.0053&0.0054&(0.5595,1.1111)&-0.0121&0.0032&(0.6891,0.9236)\\
15,15&-0.0013&0.0041&(0.5981,1.0755)&0.0122&0.0025&(0.6995,0.9168)\\
25,25&-0.0008&0.0026&(0.6454,1.0278)&-0.0056&0.0018&(0.7331,0.9076)\\
50,25&-0.0002&0.0021&(0.6657,1.0099)&-0.0040&0.0013&(0.7479,0.9018)\\
100,75&-0.0010&0.0010&(0.7211,0.9513)&-0.0005&0.0006&(0.7799,0.8850)\\
125,125&0.0002&0.0006&(0.7406,0.9343)&-0.0008&0.0005&(0.7886,0.8781)\\
150,100&0.0009&0.0007&(0.7387,0.9365)&0.0006&0.0005&(0.7888,0.8802)\\
150,150&-0.0003&0.0005&(0.7512,0.9228)&0.0003&0.0004&(0.7936,0.8760)\\
\hline
\end{longtable} 
\newpage

\begin{figure}[H]
\begin{center}
\includegraphics[height=3.5in,width=7in,angle=0]{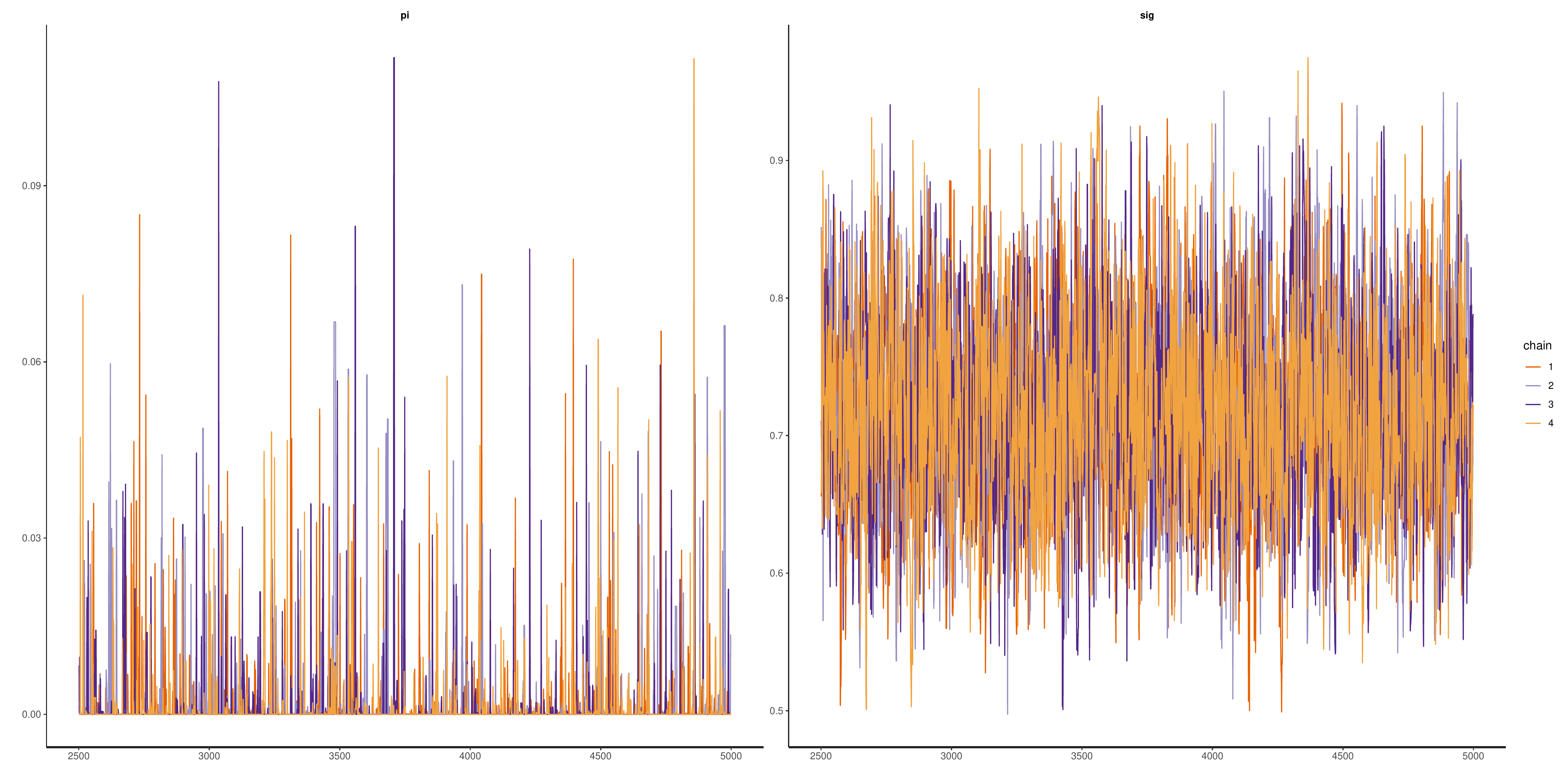}
\caption{Traceplot for $(\sigma_1,\pi_1)$ }
\end{center}
\end{figure}   
\begin{figure}[H]
\begin{center}
\includegraphics[height=3.5in,width=7in,angle=0]{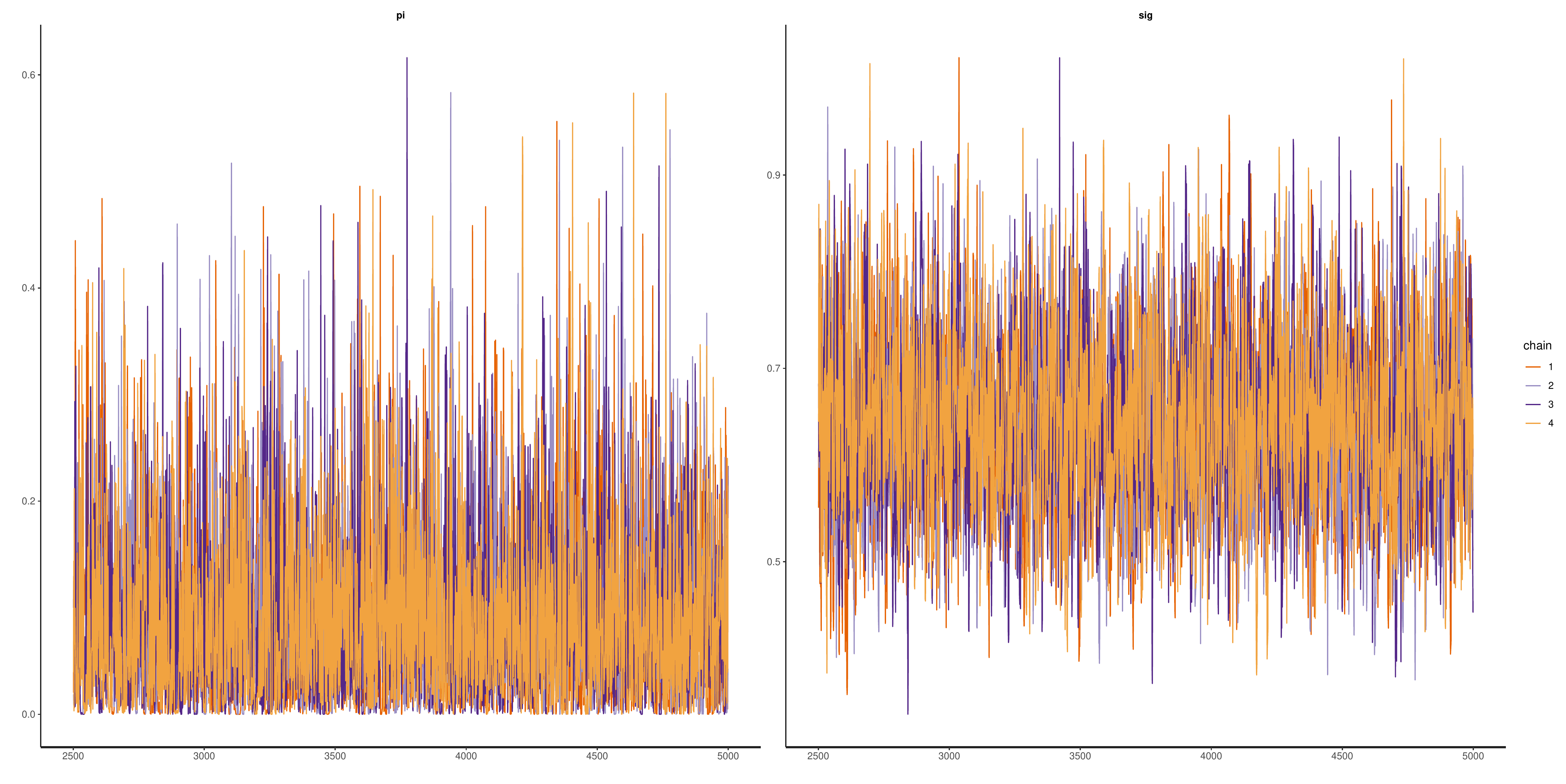}
\caption{Traceplot for $(\sigma_2,\pi_2)$ }
\end{center}
\end{figure}          
\end{document}